\documentclass[11pt]{article}
\usepackage{amssymb}
\usepackage{amsfonts}
\usepackage{amsmath}

\input{tcilatex}
\begin{document}

\begin{center}
\textbf{STABILITY OF INVERSE NODAL PROBLEM FOR ENERGY-DEPENDENT
STURM-LIOUVILLE EQUATION}

\textbf{Emrah YILMAZ\ and Hikmet KEMALOGLU}

\textbf{Firat University, Department of Mathematics, 23119, Elazig / TURKEY}

\textbf{emrah231983@gmail.com\qquad hkoyunbakan@gmail.com}
\end{center}

\begin{quote}
\textbf{Abstract: }{\footnotesize Inverse nodal problem on diffusion
operator is the problem of finding the potential functions and parameters in
the boundary conditions by using nodal data}. {\footnotesize In particular,
we solve the reconstruction and stability problems using nodal set of
eigenfunctions. Furthermore, we show that the space of all potential
functions }${\footnotesize q}${\footnotesize \ is homeomorphic to the
partition set of all asymptotically equivalent nodal sequences induced by an
equivalence relation. To show this stability which is known Lipschitz
stability, we have to construct two metric spaces and a map }$\Phi _{dif}$%
{\footnotesize \ between these spaces. We find that }$\Phi _{dif}$ 
{\footnotesize is a homeomorphism when the corresponding metrics are
magnified by the derivatives of }${\footnotesize q}${\footnotesize .
Basically, this method is similar to} {\footnotesize \cite{Tsay} and \cite%
{Cheng} which is given for Sturm-Liouville and Hill operators, respectively
and depends on the explicit asymptotic expansions of nodal points and nodal
lengths.}

\smallskip



\smallskip
\end{quote}

\textbf{\large 1. Introduction}

Inverse spectral problems is about recovering operators by using their
spectral charachteristics as spectrum, norming constants and nodal points.
Such problems have been an important research issue in mathematics and have
many applications in natural sciences. Inverse spectral problems are divided
into two parts. One of them is inverse eigenvalue problem and the other one
is inverse nodal problem. Inverse eigenvalue problems have been studied for
along time by many authors \cite{AMB}, \cite{Lev}, \cite{Lagh}, \cite{Pschl}%
, \cite{PVVRCHK}, \cite{BT}. Inverse nodal problem was first posed and
solved by McLaughlin \cite{Mc} who showed that the knowledge of a dense set
of nodal points of the eigenfunctions alone can determine the potential
function of the Sturm-Liouville problem up to a constant. Independently,
Shen and his coworkers studied the relation between the nodal points and the
density function of the string equation \cite{Shen}. Recently, many authors
have studied inverse nodal problem for different operators \cite{Hald}, \cite%
{Wang}, \cite{Yang}, \cite{YRKO}, \cite{Koyunbakan}, \cite{Koyunbakan 2}.

Consider the following boundary value problem generated by the differential
equation

\begin{equation}
Ly=y^{\prime \prime }+\left[ \lambda ^{2}-2\lambda p(x)-q(x)\right] y=0,%
\text{ }x\in (0,\pi ),  \tag{1.1}
\end{equation}%
and boundary conditions%
\begin{eqnarray}
y^{\prime }(0)-hy(0) &=&0,  \TCItag{1.2} \\
y^{\prime }(\pi )+Hy(\pi ) &=&0.  \TCItag{1.3}
\end{eqnarray}%
where $\lambda $ is spectral parameter and $h,H\in 
\mathbb{R}
;$ $q\in W_{2}^{1}\left[ 0,\pi \right] ,p\in W_{2}^{2}\left[ 0,\pi \right] $ 
\cite{Gsymv}, \cite{Yangjmp}$.$ Since $L$ is determined by its potential
functions, we identify $L$ with $L_{p,q}.$ Problem (1.1)-(1.3) are called
quadratic pencils of Schr\"{o}dinger operator (or diffusion operator). Some
versions of eigenvalue problem (1.1)-(1.3) were studied in \cite{gusein2}, 
\cite{Nab1}, \cite{Nab2}, \cite{Maksudov}, \cite{Nab3}, \cite{Yang2}, \cite%
{Pronska}, \cite{Elgiz}.

Let $\lambda _{n}$ be the $n-$th eigenvalue, $y(x,\lambda _{n})$ the
eigenfunction corresponding to the eigenvalue $\lambda _{n}$ and $%
0<x_{1}^{(n)}<x_{2}^{(n)}<...<x_{n-1}^{(n)}<\pi $ be the nodal points of the 
$n-$th eigenfunction $y(x,\lambda _{n})$. In otherwords, $%
y(x_{j}^{n},\lambda _{n})=0,j=1,2,...,n-1.$ Also let $%
I_{j}^{(n)}=[x_{j}^{(n)},x_{j+1}^{(n)}]$ be the $j-$th nodal domain of the $%
n-$th eigenvalue and let $l_{j}^{(n)}=x_{j+1}^{(n)}-x_{j}^{(n)}$ be the
nodal length. Define $x_{0}^{n}=0,x_{n}^{n}=\pi .$ We also define the
function $j_{n}(x)$ to be the largest index $j$ such that $0\leq
x_{j}^{n}\leq x.$ Thus, $j=j_{n}(x)$ if and only if $x\in \lbrack
x_{j}^{n},x_{j+1}^{n}).$ Define $X=\{x_{j}^{n}\},n\geq 0,$ $j=1,2,...n-1.$ $X
$ is called the set of nodal points of (1.1)-(1.3).

Inverse nodal problem for diffusion operator is to determine potential
functions and parameters in the boundary conditions. This type problems have
been studied by many authors \cite{koyunbakan3}, \cite{koyunbakan4}, \cite%
{buttrn}, \cite{Yangzn}, \cite{emrah}.

This study is organized as follows: in section 2, we mention some physical
and spectral properties of diffusion operator. In section 3, we give a
reconstruction formula for potential function and some important results for
the problem (1.1)-(1.3). Finally, we solve Lipschitz and high order
Lipschitz stability problems for diffusion operator in sections 4 and 5,
respectively.

\bigskip

\textbf{\large 2. Some Physical and Spectral Properties of Diffusion Equation%
}

Jaulent and Jean \cite{Jaulent} stated the actual background of diffusion
operators and discussed the inverse problem for the diffusion equation.
Also, Gasymov and Guseinov studied the spectral theory of diffusion operator 
\cite{Gsymv}.

The problem of describing the interactions between colliding particles is of
fundemental interesting in physics. In many cases, a description can be
carried out through a well known theoretical model. In particular, one is
interested in collisions of two spinless particles, and it is supposed that
the $s-$ wave binding energies and $s-$wave scattering matrix are exactly
known from collision experiments. $s-$wave Schr\"{o}dinger equation with a
radial static potential $V$ can be written as%
\begin{equation}
y^{\prime \prime }+\left[ E-V(x)\right] y=0,\text{ }x\geq 0  \tag{2.1}
\end{equation}%
where the potential function depends on energy in some way and has the
following form of energy dependence%
\begin{equation}
V(x,E)=U(x)+2\sqrt{E}Q(x).  \tag{2.2}
\end{equation}%
$U(x)$ and $Q(x)$ are complex-valued functions. (2.1) reduces to the
Klein-Gordon $s-$wave equation with the static potential $Q(x),$ for a
particle of zero mass and the energy $\sqrt{E}$ with the additional
condition $U(x)=-Q^{2}(x)$ \cite{Jaulent}$.$

The Klein Gordon equation is considered one of the most important
mathematical models in quantum field theory. The equation appears in
relativistic physics and is used to describe dispersive wave phenomena in
general. In addition, it also appears in nonlinear optics and plasma
physics. The Klein-Gordon equation arise in physics in linear and nonlinear
forms \cite{wazwaz}.

Now, we will consider Klein Gordon wave equation. After some straightforward
computations, it turns diffusion equation which is given in (1.1). This form
of Klein Gordon equation was first given in \cite{Sharma}. But here we shall
improve our understanding. Let consider following Klein Gordon wave equation%
\begin{equation}
\left[ \left( i\frac{\partial }{\partial t}-e\phi \right) ^{2}-\left( \frac{1%
}{i}\nabla -e\overrightarrow{A}\right) ^{2}\right] \Psi =m^{2}\Psi . 
\tag{2.3}
\end{equation}

This equation represents a spinless particle of charge $e$ and mass $m$ in a
scalar potential $\phi $ and vector potential $\overrightarrow{A}(r,t)$
where natural units $\hslash =c=1$ have been used. In case of $%
\overrightarrow{A}(r,t)=0$ and the scalar potential to be time independent,
equation (2.3) is reduced the following form

\begin{equation}
\left[ \nabla ^{2}-\frac{\partial ^{2}}{\partial t^{2}}-2ei\phi \frac{%
\partial }{\partial t}+e^{2}\phi ^{2}\right] \Psi =m^{2}\Psi .  \tag{2.4}
\end{equation}

If we set $\Psi (r,t)=\varphi (r)e^{-iEt},$ equation (2.4) takes the form

\begin{equation}
\left[ \nabla ^{2}+\left( E-e\phi \right) ^{2}\right] \varphi
(r)=m^{2}\varphi (r),  \tag{2.5}
\end{equation}%
where%
\begin{eqnarray*}
\Psi _{t}(r,t) &=&-\varphi (r)iEe^{-iEt} \\
\Psi _{tt}(r,t) &=&-\varphi (r)E^{2}e^{-iEt}.
\end{eqnarray*}

This equation is a well-known form of Klein-Gordon equation. The wave
equation (2.5) with a spherically symmetric potential energy may be written
in spherical coordinates

\begin{equation}
\left\{ \left[ \frac{1}{r^{2}}\frac{\partial }{\partial r}\left( r^{2}\frac{%
\partial }{\partial r}\right) +\frac{1}{r^{2}\sin \theta }\frac{\partial }{%
\partial \theta }\left( \sin \theta \frac{\partial }{\partial \theta }%
\right) +\frac{1}{r^{2}\sin ^{2}\theta }\frac{\partial ^{2}}{\partial \phi
^{2}}\right] +\left( E-e\phi \right) ^{2}\right\} \varphi =m^{2}\varphi . 
\tag{2.6}
\end{equation}

We first separate the radial and angular parts by substituting

\begin{equation*}
\varphi (r,\theta ,\phi )=R(r)Y(\theta ,\phi ),
\end{equation*}%
into equation (2.6) and dividing through by $\varphi ,$ we get

\begin{equation}
\frac{1}{R}\frac{d}{dr}\left( r^{2}\frac{dR}{dr}\right) +\left[ \left(
E-e\phi \right) ^{2}-m^{2}\right] r^{2}=-\frac{1}{Y}\left[ \frac{1}{\sin
\theta }\frac{\partial }{\partial \theta }\left( \sin \theta \frac{\partial Y%
}{\partial \theta }\right) +\frac{1}{\sin ^{2}\theta }\frac{\partial ^{2}Y}{%
\partial \phi ^{2}}\right] .  \tag{2.7}
\end{equation}

Since the left side of equation (2.7) depends only on $r$, and the right
side depends only on $\theta $ and $\phi ,$ both sides must be equal to a
constant that we call $l(l+1)$ where $l$ is orbital quantum number. Then,
the equation (2.7) gives us a radial equation with the scalar potential $%
V=e\phi $ for the Klein Gordon equation

\begin{equation}
\frac{1}{r^{2}}\frac{d}{dr}\left( r^{2}\frac{dR}{dr}\right) +\left[ \left(
E^{2}-m^{2}\right) +V(V-2E)-\frac{l(l+1)}{r^{2}}\right] R=0  \tag{2.8}
\end{equation}%
where $l=0,1,2,...$ Substituting $R=\dfrac{\psi (r)}{r}$ in equation (2.8),
we obtain

\begin{equation*}
\psi ^{\prime \prime }(r)+\left[ K^{2}+V(V-2E)-\frac{l(l+1)}{r^{2}}\right]
\psi (r)=0,
\end{equation*}%
where $E^{2}=K^{2}+m^{2}.$ In case of $l=m=0,$ we get

\begin{equation}
\psi ^{\prime \prime }(r)+\left[ V^{2}-2KV\right] \psi (r)=-K^{2}\psi (r). 
\tag{2.9}
\end{equation}

\bigskip 

There are many important spectral properties of eigenvalues and
eigenfunctions for the problem (1.1)-(1.3). We collect some of these in the
lemmas below.

\bigskip

{\large \textbf{Lemma 2.1.}} \textit{The eigenvalues of the operator }$%
L_{p,q}$ \textit{are simple. }

{\large \textbf{Proof: }}It can be easily proved by using similar way with 
\cite{Chadan}.

\bigskip

{\large \textbf{Lemma 2.2. }}\textit{The operator }$L_{p,q}$\textit{\
defined by (1.1)-(1.3) is in fact symmetric on space }$L_{2}[0,\pi ].$

{\large \textbf{Proof: }}Let $u$ and $v$ be twice differentiable functions
which satisfy the boundary conditions (1.2)-(1.3). This lemma can be esily
proved by using integration by parts on these two functions $u$ and $v$.
Considering the definition of inner product on $L_{2}[0,\pi ]$, we get%
\begin{eqnarray*}
&<&L_{p,q}\left[ u\right] ,v>=\dint\limits_{0}^{\pi }L_{p,q}\left[ u\right]
vdx=\dint\limits_{0}^{\pi }\left\{ -u^{\prime \prime }+\left[ q(x)+2\lambda
p(x)\right] u\right\} vdx \\
&=&hu(0)v(0)+v(\pi )Hu(\pi )-v(\pi )Hu(\pi )-hu(0)v(0)+\dint\limits_{0}^{\pi
}\left\{ -v^{\prime \prime }+\left[ q(x)+2\lambda p(x)\right] v\right\} udx
\\
&=&<u,L_{p,q}\left[ v\right] >.
\end{eqnarray*}

\bigskip

{\large \textbf{Lemma 2.3. }}\textit{The eigenfunctions corresponding to
different eigenvalues are orthogonal for the operator }$L_{p,q}.$

{\large \textbf{Proof: }}Let $\phi _{m}$ and $\phi _{n}$ be the
eigenfunctions corresponding to the eigenvalues $\lambda _{m}$ and $\lambda
_{n},\left( \lambda _{m}\neq \lambda _{n}\right) $ respectively$.$ If $\phi
_{k}$ denotes the $k-$th eigenfunction, we can integrate the identity 
\begin{equation*}
-\phi _{m}^{\prime \prime }\phi _{n}+\phi _{n}^{\prime \prime }\phi
_{m}+2p(x)\phi _{m}\phi _{n}\left( \lambda _{m}-\lambda _{n}\right) =\left(
\lambda _{m}^{2}-\lambda _{n}^{2}\right) \phi _{m}\phi _{n},
\end{equation*}%
to obtain 
\begin{equation*}
\dint\limits_{0}^{\pi }\left[ 2p(x)-\lambda _{m}-\lambda _{n}\right] \phi
_{m}\phi _{n}dx=0,
\end{equation*}%
from which the result follows.

\bigskip

{\large \textbf{Lemma 2.4. }}\cite{Gsymv}, \cite{Yangzn} \textit{Let }$%
\lambda _{n},n\in 
\mathbb{Z}
-\{0\}$ \textit{be the spectrum of the problem (1.1)-(1.3). It is well known
that the sequence} $\left\{ \lambda _{n}:n=1,2,...\right\} $ \textit{%
satisfies the following asymptotic formula}%
\begin{equation}
\lambda _{n}=n+c_{0}+\frac{c_{1}}{n}+\frac{c_{1,n}}{n},  \tag{2.10}
\end{equation}%
\textit{where}%
\begin{equation*}
c_{0}=\frac{1}{\pi }\dint\limits_{0}^{\pi }p(x)dx,c_{1}=\frac{1}{\pi }\left[
h+H+\frac{1}{2}\dint\limits_{0}^{\pi }\left[ q(x)+p^{2}(x)\right] dx\right]
,\dsum\limits_{n=1}^{\infty }\left\vert c_{1,n}\right\vert ^{2}<\infty .
\end{equation*}

Let consider the equation (1.1) with the initial conditions 
\begin{equation}
y(0)=0,\text{ }y^{\prime }(0)=1,  \tag{2.11}
\end{equation}%
We will denote by $\varphi (x,\lambda )$ the solution of (1.1) satisfying
the initial condition (1.2) and by $\psi (x,\lambda )$ the solution of the
same equation satisfying the initial conditions (2.11) \cite{koyunbakan3}.

\bigskip

{\large \textbf{Lemma 2.5. }}\cite{koyunbakan3} \textit{The solutions of the
problems (1.1)-(1.3) and (1.1),(1.3),(2.11) have the following forms,}%
\begin{equation}
\varphi (x,\lambda )=\cos \left( \lambda x\right) -\frac{h}{\lambda }\sin
\left( \lambda x\right) +\dint\limits_{0}^{x}\frac{\sin \left[ \lambda
\left( x-t\right) \right] }{\lambda }\left[ q(t)+2\lambda p(t)\right]
\varphi (t,\lambda )dt,  \tag{2.12}
\end{equation}%
\begin{equation}
\psi (x,\lambda )=\frac{\sin \left( \lambda x\right) }{\lambda }%
+\dint\limits_{0}^{x}\frac{\sin \left[ \lambda \left( x-t\right) \right] }{%
\lambda }\left[ q(t)+2\lambda p(t)\right] \psi (t,\lambda )dt,  \tag{2.13}
\end{equation}%
\textit{respectively.}

{\large \textbf{Lemma 2.6. }}\cite{koyunbakan3} \textit{Suppose that }$q\in
L_{1}[0,\pi ]$ \textit{and} $p\in W_{2}^{2}[0,\pi ].$ \textit{Then, the
nodal points and nodal lengths has the following asymptotic forms, as} $%
n\rightarrow \infty $

\textbf{(a)} \textit{For the problem (1.1)-(1.3),}%
\begin{equation}
x_{j}^{n}=\frac{\left( j-\dfrac{1}{2}\right) \pi }{\lambda _{n}}-\frac{h}{%
2\lambda _{n}^{2}}+\frac{1}{2\lambda _{n}^{2}}\dint\limits_{0}^{x_{j}^{n}}%
\left[ 1+\cos \left( 2\lambda _{n}t\right) \right] \left[ q(t)+2\lambda
_{n}p(t)\right] dt+o\left( \frac{1}{\lambda _{n}^{3}}\right) ,  \tag{2.14}
\end{equation}%
\begin{equation}
l_{j}^{n}=\frac{\pi }{\lambda _{n}}+\frac{1}{2\lambda _{n}^{2}}%
\dint\limits_{x_{j}^{n}}^{x_{j+1}^{n}}\left[ 1+\cos \left( 2\lambda
_{n}t\right) \right] \left[ q(t)+2\lambda _{n}p(t)\right] dt+o\left( \frac{1%
}{\lambda _{n}^{3}}\right) .  \tag{2.15}
\end{equation}

\textbf{(b)} \textit{For the problem (1.1), (1.3), (2.11),}%
\begin{equation}
x_{j}^{n}=\frac{j\pi }{\lambda _{n}}+\frac{1}{2\lambda _{n}^{2}}%
\dint\limits_{0}^{x_{j}^{n}}\left[ 1-\cos \left( 2\lambda _{n}t\right) %
\right] \left[ q(t)+2\lambda _{n}p(t)\right] dt+o\left( \frac{1}{\lambda
_{n}^{3}}\right) ,  \tag{2.16}
\end{equation}%
\begin{equation}
l_{j}^{n}=\frac{\pi }{\lambda _{n}}+\frac{1}{2\lambda _{n}^{2}}%
\dint\limits_{x_{j}^{n}}^{x_{j+1}^{n}}\left[ 1-\cos \left( 2\lambda
_{n}t\right) \right] \left[ q(t)+2\lambda _{n}p(t)\right] dt+o\left( \frac{1%
}{\lambda _{n}^{3}}\right) .  \tag{2.17}
\end{equation}

\bigskip

{\large \textbf{Lemma 2.7.}} \textit{Suppose that }$q\in L_{1}[0,\pi ].$ 
\textit{Then, for almost every} $x\in \left( 0,\pi \right) ,$ \textit{with} $%
j=j_{n}(x)$%
\begin{equation*}
\lim_{n\rightarrow \infty }\lambda
_{n}\dint\limits_{x_{j}^{n}}^{x_{j+1}^{n}}q(t)dt=q(x).
\end{equation*}

{\large \textbf{Proof: }}It can be proved by using similar way with \cite%
{Tsay}.

The above equalities are still valid even if $x_{j}^{n}$ and $x_{j+1}^{n}$
are replaced by $X_{J}^{n}$ and $X_{J+1}^{n}.$

\bigskip

\textbf{\large 3. A Reconstruction Formula for Potential Function and some
Important Results}

\bigskip

In this section, we will give a reconstruction formula for the potential
function of diffusion operator and some results. Then, we shall give
definitions of $d_{\Sigma _{dif}}$ and $d_{0}.$ We will denote $L_{k}^{n}$
(grid lengths) to go with $X_{k}^{n}$ in the same way that $l_{k}^{n}$
(nodal lengths) to go with $x_{k}^{n}.$ Also, we denote the space $\Omega
_{dif}$ as a collection of all diffusion operators $L_{p,q}$ and the space $%
\Sigma _{dif}$ as a collection of all admissible double sequences of nodes
such that corresponding functions are convergent in $L_{1}.$ A pseudometric $%
d_{\Sigma _{dif}}$ will be defined on $\Sigma _{dif}.$ Essentially, $%
d_{\Sigma _{dif}}\left( X,\overline{X}\right) $ is so close to

\begin{equation}
d_{0}\left( X,\overline{X}\right) =\overline{\lim\limits_{n\rightarrow
\infty }}n^{2}\pi \dsum\limits_{k=1}^{n-1}\left\vert L_{k}^{n}-\overline{L}%
_{k}^{n}\right\vert +n\dint\limits_{0}^{\pi }\left\vert p(x)-\overline{p}%
(x)\right\vert dx,  \tag{3.1}
\end{equation}%
where $X,\overline{X}\in \Sigma _{dif},L_{k}^{n}=X_{k+1}^{n}-X_{k}^{n}$ and $%
\overline{L}_{k}^{n}$ is defined similarly.

If we define $X\sim \overline{X}$ if and only if $d_{\Sigma _{dif}}(X,%
\overline{X})=0,$ then $\sim $ is an equivalence relation on $\Sigma _{dif}$
and $d_{\Sigma _{dif}\text{ }}$would be a metric for the partition set $%
\Sigma _{dif}^{\ast }=\Sigma _{dif}/\sim .$ Let $\Sigma _{dif1}\subset
\Sigma _{dif}$ be the subspace of all asymptotically equivalent nodal
sequences and let $\Sigma _{dif1}^{\ast }=\Sigma _{dif1}/\sim .$ Let $\Phi
_{dif}$ be the homeomorphism between the spaces $\Omega _{dif}$ and $\Sigma
_{dif1}^{\ast }.$ We call $\Phi _{dif}$ as a nodal map for diffusion
operator.

\bigskip

{\large \textbf{Theorem 3.1.}} \textit{Let }$q\in L_{1}[0,\pi ]$ \textit{and}
$p\in W_{2}^{2}[0,\pi ]$\textit{. Define} $F_{n}$ \textit{by}

\begin{equation*}
F_{n}=n\left[ \dsum\limits_{k=1}^{n-1}2n^{2}L_{k}^{n}-2n\pi -2p(x)\right] ,
\end{equation*}%
\textit{for the problem (1.1)-(1.3). Then}, $F_{n}$ \textit{converges to} $q$
\textit{poinwisely almost everywhere and also in} $L_{1}$ \textit{sense.
Furthermore, pointwise convergence holds for all the continuity points of} $%
q.$

\bigskip

{\large \textbf{Proof: }}We shall consider the asymptotic formulas for the
nodal lengths of the problem (1.1)-(1.3). Observe that by Lemma 2.6, we have%
\begin{equation*}
l_{j}^{n}\frac{\lambda _{n}}{\pi }-1=\frac{1}{2\lambda _{n}\pi }%
\dint\limits_{x_{j}^{n}}^{x_{j+1}^{n}}\left[ 1+\cos \left( 2\lambda
_{n}t\right) \right] \left[ q(t)+2\lambda _{n}p(t)\right] dt+o\left( \frac{1%
}{\lambda _{n}^{2}}\right) .
\end{equation*}

And, after some algebraic computations, we get%
\begin{equation*}
2\lambda _{n}^{2}\left( l_{j}^{n}\lambda _{n}-\pi \right) -2\lambda
_{n}p(x)=\lambda _{n}\dint\limits_{x_{j}^{n}}^{x_{j+1}^{n}}q(t)dt+o(1).
\end{equation*}

If we consider the left side of this equality with the asymptotic formula of 
$\lambda _{n},$ it yields%
\begin{eqnarray*}
2\lambda _{n}^{2}\left( l_{j}^{n}\lambda _{n}-\pi \right) -2\lambda _{n}p(x)
&=&\lambda _{n}\left( 2\lambda _{n}^{2}l_{j}^{n}-2\lambda _{n}\pi
-2p(x)\right)  \\
&=&\left[ n+o\left( \frac{1}{n}\right) \right] \left\{ 2\left[ n+o\left( 
\frac{1}{n}\right) \right] ^{2}l_{j}^{n}-2\left[ n+o\left( \frac{1}{n}%
\right) \right] \pi -2p(x)\right\}  \\
&=&n\left( 2n^{2}l_{j}^{n}-2n\pi -2p(x)\right) +o(1).
\end{eqnarray*}

Hence, to prove Theorem 3.1., it suffices to show Theorem 3.2. (b).

\bigskip

{\large \textbf{Remark: }} $F_{n}$ \textit{can be obtained similarly for the
problem (1.1),(1.3),(2.11).}

\bigskip

{\large \textbf{Theorem 3.2.}} \textit{Suppose that }$X\in \Sigma _{dif}$ 
\textit{is asymptotically nodal to }$L_{p,q}\in \Omega _{dif}.$ \textit{%
Then, we have the following.}

\textbf{(a)} \textit{For the problem (1.1)-(1.3), }%
\begin{equation*}
h=\lim_{n\rightarrow \infty }2\lambda _{n}\pi \left( j-\frac{1}{2}-\frac{%
\lambda _{n}}{\pi }X_{j}^{n}\right) .
\end{equation*}

\textbf{(b)} \textit{For almost every} $x\in \lbrack 0,\pi ]$%
\begin{equation*}
q(x)=\lim_{n\rightarrow \infty }2\lambda _{n}\left[ \lambda
_{n}^{2}l_{j}^{n}-\lambda _{n}\pi -p(x)\right] .
\end{equation*}%
where $q\in L_{1}[0,\pi ].$

{\large \textbf{Proof: }}

\textbf{(a)} From the Lemma 2.6, we know the following asymptotic formula
for the problem (1.1)-(1.3),

\begin{equation}
X_{j}^{n}=\frac{\left( j-\dfrac{1}{2}\right) \pi }{\lambda _{n}}-\frac{h}{%
2\lambda _{n}^{2}}+\frac{1}{2\lambda _{n}^{2}}\dint\limits_{0}^{X_{j}^{n}}%
\left[ 1+\cos \left( 2\lambda _{n}t\right) \right] \left[ q(t)+2\lambda
_{n}p(t)\right] dt+o\left( \frac{1}{\lambda _{n}^{3}}\right) .  \tag{3.2}
\end{equation}

After some computations in (3.2), we get%
\begin{equation*}
\lambda _{n}\pi \left( j-\frac{1}{2}-\frac{\lambda _{n}}{\pi }%
X_{j}^{n}\right) =\frac{h}{2}-\frac{1}{2}\dint\limits_{0}^{X_{j}^{n}}\left[
1+\cos \left( 2\lambda _{n}t\right) \right] \left[ q(t)+2\lambda _{n}p(t)%
\right] dt+o\left( \frac{1}{\lambda _{n}}\right) .
\end{equation*}

Since $X_{j}^{n}$ goes to zero as $n\rightarrow \infty ,$ we obtain

\begin{equation*}
h=\lim_{n\rightarrow \infty }2\lambda _{n}\pi \left( j-\frac{1}{2}-\frac{%
\lambda _{n}}{\pi }X_{j}^{n}\right) .
\end{equation*}

\textbf{(b)} Part (b) can be proved by using similar procedure with in \cite%
{koyunbakan3}.

\bigskip

{\large \textbf{Theorem 3.3.}} \textit{Given }$L_{p,q}$ \textit{in} $\Omega
_{dif},$ \textit{the set of nodal points} $\{x_{k}^{n}\}$ \textit{is also
asymptotically nodal to} $L_{p,q}$ \textit{itself}.

\bigskip

{\large \textbf{Proof: }}We can easily prove by using similar way with in 
\cite{Tsay}.

\bigskip

\textbf{\large 4. Lipschitz Stability of Inverse Nodal Problem for Diffusion
Operator}

\bigskip

In this section, we solve a Lipschitz stability problem for diffusion
operator. Lipschitz stability is about a continuity between two metric
spaces. To show this continuity, we will use a homeomorphism between these
two spaces. Stability problems were studied by many authors \cite{Tsay}, 
\cite{Cheng}, \cite{mrcnk}, \cite{mc10}, \cite{yangstab}.

\bigskip

\textbf{{\large \textbf{Definition} 4.1. }}Let $%
\mathbb{N}
^{\prime }=%
\mathbb{N}
-\{1,2\}.$

\textit{(i)}%
\begin{eqnarray*}
\Omega _{dif} &=&\{q\in L_{1}(0,\pi ):q\text{ \textit{is the potential
function of the diffusion equation}}\mathit{\}} \\
\Sigma _{dif} &=&\text{The collection of the all double sequences defined as}
\end{eqnarray*}%
\begin{equation*}
X=\left\{ X_{k}^{n}:k=1,2,...,n-2,n\in 
\mathbb{N}
^{\prime }\right\} ,0<X_{1}^{n}<X_{2}^{n}<...<X_{n-1}^{n}<\pi ,
\end{equation*}%
for each $n\in 
\mathbb{N}
${\large .}

\textit{(ii)} Let $X\in \Sigma _{dif}$ and define $X=\{X_{k}^{n}\}$ where $%
L_{k}^{n}=X_{k+1}^{n}-X_{k}^{n}$ and $I_{k}^{n}=\left(
X_{k}^{n},X_{k+1}^{n}\right) .$ We say $X$ is quasinodal to some $q\in
\Omega _{dif}$ if $X$ is an admissible sequence of nodes and satisfies (I)
and (II) below :

\textit{(I)} $X$ has the following asymptotics uniformly for $k,$ as $%
n\rightarrow \infty $%
\begin{equation*}
X_{k}^{n}=\frac{\left( k-\frac{1}{2}\right) \pi }{n}+O\left( \frac{1}{n^{2}}%
\right) ,k=1,2,...,n
\end{equation*}%
and the sequence%
\begin{equation*}
F_{n}=n\left[ \dsum\limits_{k=1}^{n-1}2n^{2}L_{k}^{n}-2n\pi -2p(x)\right] ,
\end{equation*}%
converges to $q$ in $L_{1}$ for the problem (1.1)-(1.3)$.$

\textit{(II)} $X$ has the following asymptotics uniformly for $k,$ as $%
n\rightarrow \infty $%
\begin{equation*}
X_{k}^{n}=\frac{k\pi }{n}+O\left( \frac{1}{n^{2}}\right) ,k=1,2,...,n.
\end{equation*}%
Similarly, for this case $F_{n}$ converges to $q$ in $L_{1}$ for the problem
(1.1), (1.3), (2.11).

\bigskip

\textbf{{\large \textbf{Definition} 4.2. }}Suppose that $X,\overline{X}\in
\Sigma _{dif}$ with $L_{k}^{n}$ and $\overline{L}_{k}^{n}$ as their
respective grid lengths. Let%
\begin{equation}
S_{n}\left( X,\overline{X}\right) =n^{2}\pi
\dsum\limits_{k=1}^{n-1}\left\vert L_{k}^{n}-\overline{L}_{k}^{n}\right\vert
+n\dint\limits_{0}^{\pi }\left\vert p(x)-\overline{p}(x)\right\vert dx. 
\tag{4.1}
\end{equation}

Define%
\begin{equation*}
d_{0}\left( X,\overline{X}\right) =\overline{\lim\limits_{n\rightarrow
\infty }}S_{n}\left( X,\overline{X}\right) \text{ and }d_{\Sigma
_{dif}}\left( X,\overline{X}\right) =\overline{\lim\limits_{n\rightarrow
\infty }}\frac{S_{n}\left( X,\overline{X}\right) }{1+S_{n}\left( X,\overline{%
X}\right) }.
\end{equation*}

This definition was first made by \cite{Tsay} for Sturm-Liouvile operator.
Since the function $f(x)=\dfrac{x}{1+x}$ is monotonic, we get%
\begin{equation*}
d_{\Sigma _{dif}}\left( X,\overline{X}\right) =\frac{d_{0}\left( X,\overline{%
X}\right) }{1+d_{0}\left( X,\overline{X}\right) }\in \left[ 0,\pi \right] .
\end{equation*}

Conversely, 
\begin{equation*}
d_{0}\left( X,\overline{X}\right) =\frac{d_{\Sigma _{dif}}\left( X,\overline{%
X}\right) }{1-d_{\Sigma _{dif}}\left( X,\overline{X}\right) }.
\end{equation*}

This equalities can be obtained easily.

\bigskip

{\large \textbf{Lemma 4.1. }}\textit{Let }$X,\overline{X}\in \Sigma _{dif}.$

\textbf{(a)} $d_{\Sigma _{dif}}$ \textit{is a pseudometric on} $\Sigma
_{dif}.$

\textbf{(b)} \textit{If} $X$ \textit{and} $\overline{X}$ \textit{belong to
different cases, then} $d_{\Sigma _{dif}}\left( X,\overline{X}\right) =1.$

\textbf{(c)} \textit{If} $X$ \textit{belongs to case (I) or case (II), then}%
\begin{equation*}
L_{k}^{n}=\frac{\pi }{n}+O\left( \frac{1}{n^{2}}\right) ,k=1,2,...,n.
\end{equation*}

\textbf{(d)} \textit{Let }$X$ and $\overline{X}$ \textit{belong to same case}%
.

\textbf{i)} \textit{The interval }$\delta _{n,k}$\textit{\ between the
points }$X_{k}^{n}$\textit{\ and }$\overline{X}_{k}^{n}$\textit{\ has length 
}$O(\dfrac{1}{n^{2}}),$\textit{\ as }$n\rightarrow \infty .$

\textbf{ii)} \textit{For all }$x\in (0,\pi ),$\textit{\ define }$%
J_{n}(x)=\max \{k:X_{k}^{(n)}\leq x\}\ $\textit{so that} $k=J_{n}(x)$ 
\textit{if and only if} $x\in \lbrack X_{J}^{n},X_{J+1}^{n})$. \textit{Then},%
\textit{\ }$\left\vert J_{n}(z)-\overline{J}_{n}(z)\right\vert \leq 1$%
\textit{\ for\ sufficiently large }$n$\textit{.}

\bigskip

{\large \textbf{Proof: }}It can be proved similar to \cite{Tsay}, \cite%
{Cheng}.

\bigskip

After following theorem, we can say that inverse nodal problem for diffusion
operator is Lipschitz stable.

\bigskip

{\large \textbf{Theorem 4.1.}} \textit{The metric spaces }$\left( \Omega
_{dif},\left\Vert .\right\Vert _{1}\right) $ and $\left( \Sigma _{dif1}/\sim
,d_{\Sigma _{dif}}\right) $ \textit{are homeomorphic to each other. Here} $%
\sim $ \textit{is the equivalence relation induced by }$d_{\Sigma _{dif}}.$ 
\textit{Furthermore}%
\begin{equation*}
\left\Vert q-\overline{q}\right\Vert _{1}=\frac{2d_{\Sigma _{dif}}\left( X,%
\overline{X}\right) }{1-d_{\Sigma _{dif}}\left( X,\overline{X}\right) },
\end{equation*}%
where $d_{\Sigma _{dif}}\left( X,\overline{X}\right) <1.$

\bigskip

{\large \textbf{Proof: }}In view of Lemma 4.1., we only need to consider
when $X,\overline{X}\in \Sigma _{dif}$ belong to same case. Without loss of
generality, let $X,\overline{X}$ belong to case I. We have to show%
\begin{equation*}
\left\Vert q-\overline{q}\right\Vert _{1}=2d_{0}\left( X,\overline{X}\right)
.
\end{equation*}

According to the Theorem 3.1., $F_{n}$ and $\overline{F}_{n}$ convergence to 
$q$ and $\overline{q}$, respectively. If we use definition of norm on $L_{1}$
for the potential functions, we get%
\begin{equation*}
\left\Vert q-\overline{q}\right\Vert _{1}\leq 2n^{3}\dint\limits_{0}^{\pi
}\left\vert L_{J_{n}(x)}^{n}-\overline{L}_{\overline{J}_{n}(x)}^{n}\right%
\vert dx+2n\dint\limits_{0}^{\pi }\left\vert p(x)-\overline{p}(x)\right\vert
dx+o(1)
\end{equation*}%
and after some algebraic operations

\begin{equation}
\left\Vert q-\overline{q}\right\Vert _{1}\leq 2n^{3}\dint\limits_{0}^{\pi
}\left\vert L_{J_{n}(x)}^{n}-\overline{L}_{J_{n}(x)}^{n}\right\vert
dx+2n^{3}\dint\limits_{0}^{\pi }\left\vert \overline{L}_{J_{n}(x)}^{n}-%
\overline{L}_{\overline{J}_{n}(x)}^{n}\right\vert dx+2n\dint\limits_{0}^{\pi
}\left\vert p(x)-\overline{p}(x)\right\vert dx+o(1).  \tag{4.2}
\end{equation}

Here, the integrals in the second and first terms can be written as%
\begin{equation*}
\dint\limits_{0}^{\pi }\left\vert \overline{L}_{J_{n}(x)}^{n}-\overline{L}_{%
\overline{J}_{n}(x)}^{n}\right\vert dx=o\left( \frac{1}{n^{3}}\right) ,
\end{equation*}%
and%
\begin{equation*}
\dint\limits_{0}^{\pi }\left\vert L_{J_{n}(x)}^{n}-\overline{L}%
_{J_{n}(x)}^{n}\right\vert dx=\frac{\pi }{n}\dsum\limits_{k=1}^{n-1}\left%
\vert L_{k}^{n}-\overline{L}_{k}^{n}\right\vert ,
\end{equation*}%
respectively. If we consider these equalities in (4.2), we get%
\begin{equation*}
\left\Vert q-\overline{q}\right\Vert _{1}\leq 2n^{3}o\left( \frac{1}{n^{3}}%
\right) +2n^{3}\left[ \frac{\pi }{n}\dsum\limits_{k=1}^{n-1}\left\vert
L_{k}^{n}-\overline{L}_{k}^{n}\right\vert \right] +2n\dint\limits_{0}^{\pi
}\left\vert p(x)-\overline{p}(x)\right\vert dx+o(1),
\end{equation*}%
and

\begin{equation}
\left\Vert q-\overline{q}\right\Vert _{1}\leq 2n^{2}\pi
\dsum\limits_{k=1}^{n-1}\left\vert L_{k}^{n}-\overline{L}_{k}^{n}\right\vert
+2n\dint\limits_{0}^{\pi }\left\vert p(x)-\overline{p}(x)\right\vert dx+o(1).
\tag{4.3}
\end{equation}

Similarly, we can easily obtain 
\begin{equation}
\left\Vert q-\overline{q}\right\Vert _{1}\geq 2n^{2}\pi
\dsum\limits_{k=1}^{n-1}\left\vert L_{k}^{n}-\overline{L}_{k}^{n}\right\vert
+2n\dint\limits_{0}^{\pi }\left\vert p(x)-\overline{p}(x)\right\vert dx+o(1).
\tag{4.4}
\end{equation}

Considering (4.3) and (4.4) together, it yields%
\begin{equation*}
\left\Vert q-\overline{q}\right\Vert _{1}=2n^{2}\pi
\dsum\limits_{k=1}^{n-1}\left\vert L_{k}^{n}-\overline{L}_{k}^{n}\right\vert
+2n\dint\limits_{0}^{\pi }\left\vert p(x)-\overline{p}(x)\right\vert dx.
\end{equation*}

The proof is complete after by taking limit as $n\rightarrow \infty .$

\bigskip

\textbf{\large 5. High order Lipschitz stability }

In this section, we will solve a high order stability problem for diffusion
operator. For $n\in 
\mathbb{N}
,$ 
\begin{equation*}
\Omega _{dif}(N)=\left\{ q\in L_{1}\left[ 0,\pi \right] :q\in C^{N+1}\left[
0,\pi \right] \right\} .
\end{equation*}

It has been proved in \cite{emrah} that $m-$th derivative of potential
function $q$ for Diffusion operator can be approximated by some difference
quotient of nodal length $\delta ^{m}l_{j}^{n}$ where $\delta ^{m}$ is $m-$%
th order difference quotient operator defined as \cite{Tsay}

\begin{equation*}
\delta a_{j}^{(n)}=\frac{a_{j+1}^{(n)}-a_{j}^{(n)}}{a_{j}^{(n)}}\text{ and }%
\delta ^{m}a_{j}^{(n)}=\frac{\delta ^{m-1}a_{j+1}^{(n)}-\delta
^{m-1}a_{j}^{(n)}}{a_{j}^{(n)}}.
\end{equation*}

Note that $\delta $ and $\delta ^{m}$ operators depend on the double
sequence $\{a_{i}^{(n)}\}.$

Let $\Sigma _{dif}(N)$ be the set of asymptotically equivalent nodal
sequences in $\Omega _{dif}(N).$ Let $d_{\Omega _{dif}(N)}$ and $D_{\Sigma
_{dif}(N)}$ be some metrics on $\Omega _{dif}(N)$ and $\Sigma _{dif}^{\ast
}(N)$ respectively magnified by the $L_{1}$ norms of derivatives of the
potential functions. We find that the nodal map $\Phi _{dif}$ is still a
homeomorphism under these strengthed metrics.

If we define $X\sim _{N}\overline{X}$ if and only if $D_{\Sigma
_{dif}(N)}\left( X,\overline{X}\right) =0,$ then $\sim _{N}$ is an
equivalence relation on $\Sigma _{dif}(N).$ Hence $d_{\Omega _{dif}(N)}$ is
a metric on $\Omega _{dif}(N).$

\bigskip

\textbf{{\large \textbf{Definition} 5.1. }} Let $X,\overline{X}\in \Sigma
_{dif}$ and $X_{k}^{n}=\overline{X}_{k}^{n}=\pi $ for $k>n$. For $%
m=1,2,...,N $, let%
\begin{eqnarray}
S_{m,n}\left( X,\overline{X}\right) &=&\lambda _{n}^{\frac{1}{2}%
}\dsum\limits_{k=1}^{n-m-2}\left\vert \delta ^{m}L_{k}^{n}-\delta ^{m}%
\overline{L}_{k}^{n}\right\vert +\lambda
_{n}\dint\limits_{X_{1}^{n}}^{X_{n-m-1}^{n}}\left\vert \delta ^{m}\overline{p%
}(x_{\overline{J}})-\delta ^{m}p(x_{J})\right\vert dx  \notag \\
&&+\lambda _{n}\dint\limits_{X_{1}^{n}}^{X_{n-m-1}^{n}}\left\vert \overline{p%
}^{(m)}(x)-p^{(m)}(x)\right\vert dx.  \TCItag{5.1}
\end{eqnarray}

Define%
\begin{equation*}
d_{m}\left( X,\overline{X}\right) =\overline{\lim\limits_{n\rightarrow
\infty }}S_{m,n}\left( X,\overline{X}\right) \text{ and }d_{\Sigma
_{dif}(m)}\left( X,\overline{X}\right) =\overline{\lim\limits_{n\rightarrow
\infty }}\frac{S_{m,n}\left( X,\overline{X}\right) }{1+S_{m,n}\left( X,%
\overline{X}\right) }.
\end{equation*}

We have to use the following reconstruction formula for the function $%
q^{(m)} $ in Theorem 5.1. to prove Theorem 5.3.

\bigskip

{\large \textbf{Theorem 5.1.}} \cite{emrah} \textit{Let }$q$ \textit{is real
valued}, $(n+1)-$\textit{th order continuous function from the class} $%
L_{1}[0,\pi ]$ \textit{for} $N\geq 1$ \textit{in (1.1)}$,$ \textit{and let} $%
j=j_{n}(x)$ \textit{for each} $x\in \lbrack 0,\pi ].$ \textit{Then, as} $%
n\rightarrow \infty $%
\begin{equation*}
q(x)=\lambda _{n}\left[ 2\lambda _{n}^{2}l_{j}^{n}-2\pi \lambda _{n}-2p(x)%
\right] +O\left( \frac{1}{n}\right) ,
\end{equation*}%
\textit{and, for all} $m=1,2,...,N$%
\begin{equation*}
q^{(m)}(x)=\frac{2\lambda _{n}^{3/2}}{\pi }\delta ^{m}l_{j}^{n}-2\lambda
_{n}\delta ^{m}p(x_{j})-2\lambda _{n}p^{(m)}(x)+O(1).
\end{equation*}

{\large \textbf{Theorem 5.2. }}The metric spaces $(\Omega
_{dif}(N),d_{\Omega _{dif}(N)})$ and $(\Sigma _{dif1}(N)/\sim _{N},D_{\Sigma
_{dif}(N)})$ are homeomorphic to each other, where $\sim _{N}$ is the
equivalence relation induced by $D_{\Sigma _{dif}(N)}.$

\bigskip

\textbf{\large Proof: }In view of Lemma 4.1. we only need to consider $X,%
\overline{X}\in \Sigma _{dif}$ belong to same case. We shall show that $%
\left\Vert q^{(m)}-\overline{q}^{(m)}\right\Vert =2d_{m}(X,\overline{X}).$
Hence, to prove Theorem 5.2., it suffices to prove Theorem 5.3.

\bigskip

After Theorem 5.3., we can say that the inverse nodal problem for the
diffusion operator is high order Lipschitz stable.

\bigskip

{\large \textbf{Theorem 5.3.}}\textit{\ Suppose that }$q$\textit{\ and }$%
\overline{q}$\textit{\ are both }$N$\textit{th order continuous functions
from the class }$L_{1}[0,\pi ];$ $L_{p,q},\overline{L}_{\overline{p},%
\overline{q}}$ \textit{belong to same case.} \textit{\ Let }$X$\textit{\ and 
}$\overline{X}$\textit{\ be the corresponding asymptotically equivalent
nodal sequences. Then, for all }$m=1,2,...,N$%
\begin{equation*}
\left\Vert q^{(m)}-\overline{q}^{\text{ }(m)}\right\Vert =2d_{m}(X,\overline{%
X}).
\end{equation*}

{\large \textbf{Proof: }} By Theorem 5.1 and Lemma 2.6, we can write 
\begin{equation*}
q^{(m)}-\overline{q}^{(m)}=\frac{2\lambda _{n}^{3/2}}{\pi }\left[ \delta
^{m}L_{J}^{n}-\delta ^{m}\overline{L}_{\overline{J}}^{n}\right] +2\lambda
_{n}\left[ \delta ^{m}\overline{p}(x_{\overline{J}})-\delta ^{m}p(x_{J})%
\right] +2\lambda _{n}\left[ \overline{p}^{(m)}(x)-p^{(m)}(x)\right]
+o\left( 1\right) ,
\end{equation*}

and%
\begin{equation*}
\left\vert q^{(m)}-\overline{q}^{(m)}\right\vert =\left\vert \frac{2\lambda
_{n}^{3/2}}{\pi }\left[ \delta ^{m}L_{J}^{n}-\delta ^{m}\overline{L}_{%
\overline{J}}^{n}\right] +2\lambda _{n}\left[ \delta ^{m}\overline{p}(x_{%
\overline{J}})-\delta ^{m}p(x_{J})\right] +2\lambda _{n}\left[ \overline{p}%
^{(m)}(x)-p^{(m)}(x)\right] +o(1)\right\vert .
\end{equation*}

Then, by using the definition of norm on $L_{1}$ space, we obtain%
\begin{equation*}
\left\Vert q^{(m)}-\overline{q}^{(m)}\right\Vert
_{1}=\dint\limits_{X_{1}^{n}}^{X_{n-m-1}^{n}}\left\vert q^{(m)}(x)-\overline{%
q}^{(m)}(x)\right\vert dx,
\end{equation*}

and%
\begin{eqnarray}
\left\Vert q^{(m)}-\overline{q}^{(m)}\right\Vert _{1}
&=&\dint\limits_{X_{1}^{n}}^{X_{n-m-1}^{n}}\left\vert \frac{2\lambda
_{n}^{3/2}}{\pi }\left[ \delta ^{m}L_{J}^{n}-\delta ^{m}\overline{L}_{%
\overline{J}}^{n}\right] +2\lambda _{n}\left[ \delta ^{m}\overline{p}(x_{%
\overline{J}})-\delta ^{m}p(x_{J})\right] \right.  \notag \\
&&\left. +2\lambda _{n}\left[ \overline{p}^{(m)}(x)-p^{(m)}(x)\right]
+o(1)\right\vert dx  \TCItag{5.2}
\end{eqnarray}%
when $n$ is sufficiently large. By using the property of triangle inequality
in (5.2), we get%
\begin{eqnarray*}
\left\Vert q^{(m)}-\overline{q}^{(m)}\right\Vert _{1} &\leq &\frac{2\lambda
_{n}^{3/2}}{\pi }\dint\limits_{X_{1}^{n}}^{X_{n-m-1}^{n}}\left\vert \delta
^{m}L_{J}^{n}-\delta ^{m}\overline{L}_{\overline{J}}^{n}\right\vert
dx+2\lambda _{n}\dint\limits_{X_{1}^{n}}^{X_{n-m-1}^{n}}\left\vert \delta
^{m}\overline{p}(x_{\overline{J}})-\delta ^{m}p(x_{J})\right\vert dx \\
&&+2\lambda _{n}\dint\limits_{X_{1}^{n}}^{X_{n-m-1}^{n}}\left\vert \overline{%
p}^{(m)}(x)-p^{(m)}(x)\right\vert dx+o(1),
\end{eqnarray*}%
and after some computations, 
\begin{eqnarray*}
\left\Vert q^{(m)}-\overline{q}^{(m)}\right\Vert _{1} &\leq &\frac{2\lambda
_{n}^{3/2}}{\pi }\dint\limits_{X_{1}^{n}}^{X_{n-m-1}^{n}}\left\vert \delta
^{m}L_{J}^{n}-\delta ^{m}\overline{L}_{J}^{n}\right\vert dx+\frac{2\lambda
_{n}^{3/2}}{\pi }\dint\limits_{X_{1}^{n}}^{X_{n-m-1}^{n}}\left\vert \delta
^{m}\overline{L}_{J}^{n}-\delta ^{m}\overline{L}_{\overline{J}%
}^{n}\right\vert dx \\
&&+2\lambda _{n}\dint\limits_{X_{1}^{n}}^{X_{n-m-1}^{n}}\delta ^{m}\overline{%
p}(x_{\overline{J}})-\delta ^{m}p(x_{J})dx+2\lambda
_{n}\dint\limits_{X_{1}^{n}}^{X_{n-m-1}^{n}}\left\vert \overline{p}%
^{(m)}(x)-p^{(m)}(x)\right\vert dx+o(1).
\end{eqnarray*}

\bigskip

Now, by using [\cite{Law}, Lemma 2.1] and Lemma 4.1., we obtain

\begin{equation*}
\dint\limits_{X_{1}^{n}}^{X_{n-m-1}^{n}}\left\vert \delta
^{m}L_{J}^{n}-\delta ^{m}\overline{L}_{J}^{n}\right\vert
dx=\dsum\limits_{k=1}^{n-m-2}\left\vert \delta ^{m}L_{k}^{n}-\delta ^{m}%
\overline{L}_{k}^{n}\right\vert \left\vert L_{k}^{n}\right\vert ,
\end{equation*}%
and%
\begin{equation*}
\dint\limits_{X_{1}^{n}}^{X_{n-m-1}^{n}}\left\vert \delta ^{m}\overline{L}%
_{J}^{n}-\delta ^{m}\overline{L}_{\overline{J}}^{n}\right\vert dx=O\left( 
\frac{1}{n^{4}}\right) .
\end{equation*}

\bigskip

Moreover, $L_{k}^{n}=\dfrac{\pi }{n}+O\left( \frac{1}{n^{2}}\right) .$ Thus,%
\begin{eqnarray*}
\frac{2\lambda _{n}^{3/2}}{\pi }\dint\limits_{X_{1}^{n}}^{X_{n-m-1}^{n}}%
\left\vert \delta ^{m}L_{J}^{n}-\delta ^{m}\overline{L}_{\overline{J}%
}^{n}\right\vert dx &\leq &\frac{2\lambda _{n}^{3/2}}{\pi }O\left( \frac{1}{%
n^{4}}\right) +\frac{2\lambda _{n}^{3/2}}{\pi }\dsum\limits_{k=1}^{n-m-2}%
\left\vert \delta ^{m}L_{k}^{n}-\delta ^{m}\overline{L}_{k}^{n}\right\vert
\left\vert \dfrac{\pi }{n}+O\left( \frac{1}{n^{2}}\right) \right\vert \\
&=&o\left( \frac{1}{n^{5/2}}\right)
+2n^{1/2}\dsum\limits_{k=1}^{n-m-2}\left\vert \delta ^{m}L_{k}^{n}-\delta
^{m}\overline{L}_{k}^{n}\right\vert +\frac{2}{\pi }n^{3/2}nO(n^{-2})o\left( 
\frac{1}{n^{3}}\right) \\
&=&2n^{1/2}\dsum\limits_{k=1}^{n-m-2}\left\vert \delta ^{m}L_{k}^{n}-\delta
^{m}\overline{L}_{k}^{n}\right\vert +o\left( \frac{1}{n^{5/2}}\right) .
\end{eqnarray*}

Therefore, 
\begin{eqnarray}
\left\Vert q^{(m)}-\overline{q}^{(m)}\right\Vert _{1} &\leq
&2n^{1/2}\dsum\limits_{k=1}^{n-m-2}\left\vert \delta ^{m}L_{k}^{n}-\delta
^{m}\overline{L}_{k}^{n}\right\vert +2\lambda
_{n}\dint\limits_{X_{1}^{n}}^{X_{n-m-1}^{n}}\left\vert \delta ^{m}\overline{p%
}(x_{\overline{J}})-\delta ^{m}p(x_{J})\right\vert dx  \notag \\
&&+2\lambda _{n}\dint\limits_{X_{1}^{n}}^{X_{n-m-1}^{n}}\left\vert \overline{%
p}^{(m)}(x)-p^{(m)}(x)\right\vert dx+o\left( 1\right) .  \TCItag{5.3}
\end{eqnarray}

If we take limit as $n\rightarrow \infty $ and use the metric definition, we
obtain%
\begin{equation}
\left\Vert q^{(m)}-\overline{q}^{(m)}\right\Vert \leq 2d_{m}(X,\overline{X}).
\tag{5.4}
\end{equation}

Using similar procedure, 
\begin{equation}
\left\Vert q^{(m)}-\overline{q}^{(m)}\right\Vert \geq 2d_{m}(X,\overline{X}).
\tag{5.5}
\end{equation}

Then, considering (5.4) and (5.5) together, we get%
\begin{equation*}
\left\Vert q^{(m)}-\overline{q}^{(m)}\right\Vert =2d_{m}(X,\overline{X}).
\end{equation*}

This completes the proof.

\end{document}